\documentclass[12pt]{amsart}\usepackage{amssymb}

\title[p-adic local dynamics] {A p-adic approach to local analytic dynamics: analytic flows and analytic maps tangent to the identity}

\subjclass[2000]{Primary 30D05, 32P05}
\keywords{conjugacy, normal form, non-archimedean, p-adic, formal, analytic, holomorphic}

\author{ Adrian Jenkins \and Steven Spallone }
\address{Department of Mathematics, Purdue University, West Lafayette,
IN, 47906} \email{majenkin@math.purdue.edu}
\address{Department of Mathematics, Purdue University, West Lafayette,
IN, 47906} \email{sspallon@math.purdue.edu}

\newtheorem{thm}{Theorem}[section]
\newtheorem{cor}[thm]{Corollary}
\newtheorem{prop}[thm]{Proposition}
\newtheorem{define}[thm]{Definition}

\newtheorem{lemma}[thm]{Lemma}

\newcommand{\nc}{\newcommand}
\nc{\dmo}{\DeclareMathOperator}
\nc{\ul}{\underline}
\nc{\N}{\mathbb{N}}
\nc{\Q}{\mathbb{Q}}
\nc{\Z}{\mathbf{Z}}
\nc{\R}{\mathbf{R}}
\nc{\C}{\mathbf C}
\nc{\OK}{\mathcal{O}_0^K}
\nc{\eps}{\varepsilon}
\nc{\beq}{\begin{equation}}
\nc{\ra}{\rightarrow}
\nc{\gm}{\gamma}
\nc{\ol}{\overline}

\dmo{\ord}{ord}

\begin{document}

\begin{abstract}
In this note, we will consider the question of local equivalence of analytic functions which fix the origin and are tangent to the identity, as well as the question of flows of analytic vector fields. All mappings and equivalences are considered in the non-archimedean context e.g. all norms can be considered $p$-adic norms. We show that any two mappings $f$ and $g$ which are formally equivalent are also analytically equivalent, and we show that analytic vector fields generate analytic flows. We consider the related questions of roots and centralizers for analytic mappings. In this setting, anything which can be done formally can also be done analytically.
\end{abstract}

\maketitle

\section{Introduction}\setcounter{equation}{0}\label{section:s1}
The goal of this paper is to consider the local analytic equivalence of mappings $f$ which are tangent to the identity, but whose convergence is with respect to a non-archimedean norm $|\cdot |$ (for example, $p$-adic norms), as well as the generation of analytic flows for analytic vector fields. We discover the interesting fact that the non-archimedean case yields very simple invariants for the analytic classification, as opposed to the substantially more difficult (both in construction and interpretation) invariants present in the theory for ${\mathbf C}$.

The question of the local conjugacy classes of holomorphic mappings (analytic in $\C$) has a long history. In 1884, the first such results regarding equivalences were proven by Koenigs \cite{koenigs}. In particular, it was shown that given an holomorphic function $f(z)=az+O(z^{2})$ defined in a neighborhood $U$ of ${\mathbf C}$, where $|a|\neq 1$, then in a neighborhood $V\subseteq U$ of the origin, there is a conformal mapping $h(z)=z+O(z^{2})$ such that $(h\circ f\circ h^{-1})(z)=az$. Thus, any such mapping can be linearized, and moreover, the linearizing biholomorphism $h$ can be taken to be tangent to the identity. 

Of course, one cannot naively expect linearization if $|a|=1$; as an obvious (and important) example, if $a=1$, then linearization is   impossible for any function $f(z)\neq z$. Thus, an interesting problem is to determine the invariants present in such a classification of mappings tangent to the identity. It is easy to acquire formal invariants for this equivalence (we review this in Section \ref{section:s2}). Nonetheless, the holomorphic classification remains very delicate; after initial attempts by Fatou \cite{Fatou} in the late 1910's to determine the invariants, the problem remained unsolved until the early 1980's, when Ecalle \cite{Ecalle} and Voronin \cite{Voronin} independently developed the moduli space of invariants for such mappings (see also the work of Malgrange \cite{Malgrange} and Il'yashenko \cite{Ilya} for a different approach). We note here that such holomorphic classification relied on a topological conjugacy, provided independently by Camacho \cite{Camacho} and Shcherbakov \cite{Shch}. 

One may also address the related question of analytic flows for analytic vector fields. In particular, it is a desirable property that a mapping tangent to the identity may be embedded in the flow of an analytic vector field; for example, such mappings admit root extraction of any order, and have an easily identifiable centralizer. In ${\mathbf C}$, this can be done both formally and topologically (or even smoothly away from the origin; see \cite{Jenkins2}), but holomorphically, this is impossible. In fact, by the realization theorems of Ecalle and Voronin mentioned above, it is incredibly rare for a holomorphic germ tangent to the identity to be the time-$t$ map of a holomorphic vector field. Despite this theory, however, it is quite difficult to determine if a particular mapping $f$ is embeddable.

Strangely enough, the following problem is still very much open: given two mappings $f$ and $g$, are these two mappings equivalent via an analytic map which is tangent to the identity? While the analytic classifications cited above provide the theoretical invariants present, as Ahern and Rosay note \cite{A-R}, they are too difficult to be computed for even the most trivial of mappings. There are partial results present (for example, it is known which {\it entire} functions are analytically equivalent to $f(z)=\frac{z}{1-z}$, and which are equivalent to $f(z)=z+z^{2}$), but as a whole, the problem is poorly understood. It is easy to construct formal power series $H$ conjugating $f$ to $g$, but showing that any such power series converges (or diverges) is generally very difficult. 

If we restrict ourselves to the case where the mappings $f$ and $g$ have {\it rational} coefficients, then we may often take the conjugating power series $H$ to have rational coefficients.  In this situation it is natural to study the {\it $p$-adic} convergence of $H$ for a given prime $p$.  Roughly speaking, this analyzes the largest power of $p$ which divides the denominators of the coefficients of $H$.  We may view the rational coefficients of $H$ as sitting inside the $p$-adic completion ${\mathbb Q}_p$ of $\Q$ rather than the archimedean completion $\R$, and do our work there.  This study complements the classical question of holomorphic convergence; it is an instance of what is popularly known as the Lefschetz principle, which roughly says that interesting questions for real and complex numbers should have interesting analogues in the $p$-adic setting.  This principle has found application in harmonic analysis, algebraic number theory, and more recently in dynamical systems (see for example  \cite{Benedetto}, \cite{HC}, etc).

Indeed, for power series with coefficients in a complete, non-archimedean valued field, it becomes reasonable to test for the convergence of a given conjugating map $H$. The reason for this is two-fold: first, a series $\sum a_{n}$ converges with respect to a non-archimedean norm if and only if $a_{n}\rightarrow 0$ as $n\rightarrow \infty$. Secondly (in a sense to be made precise later), the convergence of a power series depends solely on the decay of denominators - growth in the numerator is not detrimental to convergence.

For any field $K$ with norm $|\cdot |$, denote the ring of absolutely convergent power series centered at $0$ with coefficients in $K$ as $\OK$. In this paper, we consider $K$ of characteristic $0$ which are complete, non-archimedean valued fields, but our interest will be in the field of $p$-adic numbers ${\mathbb Q}_{p}$, any finite extension of ${\mathbb Q}_{p}$, and the analytic completion of any infinite algebraic extension of ${\mathbb Q}_{p}$.

The main results of this paper are the following:

\begin{thm}\label{thm1.1}
Fix any prime $p$, and let $f\in \OK$ be an analytic function which is tangent to the identity, $f(x)=x+a_{m}x^{m}+\cdots $. Write $\widetilde K=K[\sqrt[m-1]{a_m}]$. Then, there is a $\mu \in \widetilde K$ and an analytic function $h\in {\mathcal O}_{0}^{\widetilde K}$ tangent to the identity so that $h\circ f\circ h^{-1}(x)=x+x^{m}+\mu x^{2m-1}$. Moreover, $m$ and $\mu $ are analytic invariants for $f$.
\end{thm}

(The algebraic technicality of adjoining a root is unfortunately necessary, although not if $K$ is algebraically closed. We will usually drop the tilde in practice - this should cause no confusion). 

In other words, the formal and analytic classifications agree in the non-archimedean setting. This is in stark contrast to the analytic classification in ${\mathbf C}$. It is worth noting that, since analytic functions with respect to a non-archimedean norm are continuous, this shows that formal equivalence does indeed imply topological equivalence, which is consistent with the theory in ${\mathbf C}$. 

We also address the issue of analytic flows, and the next theorem shows  that the existence of a formal flow implies the existence of an analytic one.

\begin{thm}\label{thm1.2}
Consider the analytic vector field $V(x)=\left( \sum_{n=1}^{\infty} b_{n}x^{n} \right) \frac{\partial}{\partial x}  $ with coefficients in a non-archimedean field $K$, and suppose that this field generates a formal flow $\{ T^{t}_{V}(x)\} $. Then, this flow is also analytic.
\end{thm}

Theorem \ref{thm1.2} may be known, but it is included with a complete proof, since it can be used together with Theorem \ref{thm1.1} in order to solve the problem of centralizers and root extraction for mappings tangent to the identity. 

The structure of the paper is as follows: Section \ref{section:s2} discusses basic results and notation for non-archimedean analysis and local dynamics. Section \ref{section:s3} proves Theorem \ref{thm1.1} in the case where $f$ has integral coefficients.
This is not only of independent interest, but also indicates some steps necessary for the proof of Theorem \ref{thm1.1}.  We postpone the proof of Theorem \ref{thm1.1} until Section \ref{section:s4}; instead Section \ref{section:s2.5} 
is used to give a proof of Theorem \ref{thm1.2}.

Finally, Section \ref{section:s4} is devoted to a complete proof of the Theorem \ref{thm1.1}, together with results on centralizers and root extractions.

Theorems \ref{thm1.1} and \ref{thm1.2} provide a complete analytic classification of mappings/flows which are tangent to the identity and convergent with respect to a non-archimedean norm. As mentioned, this shows that formally equivalent mappings are also topologically equivalent. In a future work, the authors plan to give a complete topological classification of analytic mappings tangent to the identity, and to study what smoothness conditions may be imposed on such a conjugating map. Recently, Jenkins \cite{Jenkins1} has given a full formal classification of so-called semi-hyperbolic mappings in ${\mathbf C}^{n}$, and has shown that the formal classification differs wildly from the holomorphic one. There are few field restrictions on the formal classification; the techniques used there would work, for the most part, if the coefficients lay in any field of characteristic $0$. It would be of interest to determine whether the two classifications agree if one considers non-archimedean norms $|\cdot |$. Finally, we do not consider the case of fields with characteristic $p$; the methods used here will fail in that setting. 

This work was written while both authors were Research Assistant Professors at Purdue University. The authors are thankful for the support of the department.

\section{Preliminaries}\setcounter{equation}{0}\label{section:s2}
This section is devoted to an explanation of the non-archimedean setting in which we work, together with some basic notions of formal dynamics. We also take the opportunity to fix some notation. In the introduction and throughout the paper, we have used the convention that if the norm is archimedean, then we will write sets in the bold style (i.e. $ {\mathbf R}, {\mathbf C}$, etc.), whereas if the norm is non-archimedean, we will write in blackboard bold style (i.e. ${\mathbb Q}, {\mathbb Q}_{p}, {\mathbb C}_{p}$, etc.). Power series with coefficients in ${\mathbf C}$ will be denoted $f(z)$ as usual, while power series with coefficients in a non-archimedean field $K$ will be denoted $f(x)$. 

\subsection{Non-archimedean Fields}
The bulk of this paper involves not the complex numbers $\C$ but rather a non-archimedean complete (nontrivial) valued field $K$ of characteristic $0$. We give a survey of the pertinent facts.  For proofs see \cite{Schikhof} or \cite{Serre}. 

\begin{define} Let $K$ be a field.  A non-archimedean valuation (or  norm) on $K$ is a map $| \cdot |: K \rightarrow \R$ satisfying the following rules, for all $x,y \in K$:
\begin{enumerate}
\item $|x| \geq 0, |x|=0$ if and only if $x=0$.
\item $|x+y| \leq \max \{|x|,|y| \}$.
\item $|xy|=|x||y|$.
\end{enumerate}
The pair $(K,|\cdot |)$ is a non-archimedean valued field.
\end{define}

We will simply write $K$ when the valuation is implicit.
Of course the usual absolute value in $\C$ does not satisfy the second condition.
The constant valuation, $|x|=1$ for all $x\neq 0$, is called trivial.  We do not consider these.

Let $K=\Q$ and choose a prime $p \in \Z$ and a real number $0<\alpha <1$.  Consider the map 
\beq
\left| \frac{m}{n}\right| _{p,\alpha}=\alpha^{\ord_p(m)-\ord_p(n)},
\end{equation}
where $\ord_p(n)$ is the exponent of $p$ in the prime factorization of $n$.
Then $|\cdot |_{p,\alpha}$ is a non-archimedean valuation on $\Q$.

The following is a well-known theorem of Ostrowski:

\begin{prop} Any nontrivial non-archimedean valuation on $\Q$ is of the form $| \cdot|_{p,\alpha}$ for some $p$ and $\alpha$ as above. \end{prop}

Given a valuation on a field $K$, there is a natural topology on $K$ compatible with $|\cdot |$.
We define it in the usual way with balls.
\begin{define} Given a positive number $r \in \R$, and $x \in K$, define $B_r(x)=\{ y \in K : |x-y| \leq r \}$. \end{define}

Then we give $K$ the topology generated by the basis $\{ B_r(x) :r \in \R,x \in K \}$. For a given $p$, the topology  of $({\mathbb Q}, |\cdot |_{p,\alpha })$ does not depend on the choice of $\alpha $.

\begin{define} $\Delta =B_1(0).$ \end{define}

Note that $\Delta $ is a subring of $K$ by the definition of valuation; we will refer to $\Delta $ as the ring of integers of $K$.  

A non-archimedean valued field $K$ is considered complete if it is complete as a topological space. Recall that if $\widetilde K$ is a finite-degree field extension of a complete, non-archimedean valued field $K$, then the norm $|\cdot |$ on $K$ extends uniquely to $\widetilde K$, and furthermore, this extension is complete. In particular, if $\alpha $ is any algebraic element over such a field $K$, then $K[\alpha ]$ is a complete, non-archimedean valued field.

For example, $(\Q, |\cdot |_{p,\alpha})$ is not complete, being countable.  In fact the completion of any such $K$ will be a complete non-archimedean valued field.  The completion of $(\Q, |\cdot |_{p,1/p})$ is called $\Q_p$.  Note that $|p|=\frac{1}{p}$ in this case.

From now on we take $K$ to be a non-archimedean complete valued field with characteristic $0$.  In this case $\Q$ is a subfield, and becomes a valued field by restriction of $| \cdot |$.

We record a simple lower estimate on $|n!|$ in this context. 

\begin{prop}\label{factorial estimate}
If the valuation of $K$ restricts trivially to $\Q$ then $|n!|=1$.  Otherwise, $|n!|=|n!|_{p,\alpha} \geq \alpha^{n}$. 
\end{prop}

\begin{proof} The only thing to prove is the last inequality.  It is well-known that $\ord_p(n!)=\frac{n-S_n}{p-1}$, where $S_n$ is the sum of the digits of $n$ in base $p$.  Therefore $\ord_p(n!) \leq n$, and the result follows. \end{proof}

Since the valuation on $K$ is nontrivial, there is an element $\pi \in K$ with $0<|\pi|<1$.  Since $\R$ {\it is} an archimedean field, 
for every $\eps >0$ there is a $k \in {\mathbf N}$ so that if $q=\pi^k$, then $ |q|< \eps$.

Finally, we would like to point out that any algebraically closed field of characteristic $0$ with the same cardinality as $\R$ is isomorphic as a field to $\C$, by transcendence theory.  This applies, for example, to the algebraic closure $\ol{\Q}_p$ and its completion $\mathbb{C}_p$.  This means that much of the formal algebraic theory of $\C$ applies to a general nonarchimedean field $K$ of characteristic $0$. 

Of course there is no reason to expect any topological relationship.

\subsection{Power Series} We denote the ring of formal power series $K[[x]]$ as usual.

An interesting feature of non-archimedean analysis is the following: a series $\sum_n a_n$ converges if and only if $a_n \to 0$.
  
Given a power series $f(x)=\sum_n a_nx^n \in K[[x]]$, its radius of convergence about $0$ is given by 
\begin{equation}\label{radius of convergence}
\rho=\left(\limsup_{n \ra \infty} \sqrt[n]{ |a_n|} \right)^{-1}
\end{equation}

\begin{define} The power series $f(x) \in K[[x]]$ is called locally analytic at $0$ if $\rho>0$. 
The set of such functions is denoted $\OK$.
\end{define}

For example if $\gm \in K$ with $|\gm|=c$, and $a_n=\gm^n$, then $\rho=\frac{1}{c}$.
On the other hand, if $c >1$ and $a_n=\gm^{n^2}$, then $\rho=0$.
Therefore if $K=\Q_p$ the power series 
\beq
f(x)=\sum_n \frac{1}{p^{n^2}}x^n
\end{equation}
is not in $\OK$. 

As usual, if the linear term $a_1 \neq 0$, then $f$ will be formally invertible, in the sense that there is a unique power series $g(x) \in K[[x]]$ with $(f \circ g)(x)=(g \circ f)(x)=x$.  In particular, power series of the form $f(x)=x+O(x^{2}) $ are invertible. Moreover,  an implicit function theorem  implies that if $f$ is locally analytic, then the formal inverse $g$  is itself locally analytic.

Let $f(x)=x+a_{m}x^{m}+O(x^{m+1})$ be a power series in $K[|x|]$, where $K$ is any field of characteristic $0$. We will write $f\circ g$ to be the composition of $f$ and $g$, while writing $fg$ to mean the standard multiplicative, pointwise product. Furthermore, given $n\in {\mathbf Z}$, we write $f^{\circ n}$ to be the $n$th iterate of $f$, and write $f^{n}$ to be the $n$th multiplicative power of $f$. Given two such power series $f$ and $g$, we say that $f$ and $g$ are equivalent (or conjugate) if there is an $h$ satisfying $h\circ f\circ h^{-1}=g$. We are deliberately vague here - as mentioned in the introduction, the degree of smoothness on the map $h$ will have a huge effect on the equivalence classes present. In this paper, we will concern ourselves with two cases: $h$ can be a formal power series, or an analytic one (if $h$ is analytic, then obviously $K$ will have some associated norm).  

By considering the conjugating map $x\mapsto (\root{m-1}\of{a_{m}})x$, we  assume that $a_{m}=1$, and this assumption will be present throughout the paper. We show here that any such mapping may be reduced formally to $f_{0,m}(x)=x+x^{m}+\mu x^{2m-1}$, and so the numbers $m$ and $\mu $ provide formal invariants for the mapping $f$. The proof of this fact is known to many, and is impossible to ascribe to a single source. We include the proof here, however, as the mapping constructed will always converge in the non-archimedean setting (as we shall show later).

\begin{prop}\label{formal equivalence}
Let $f\in K[|x|]$ have the form
\begin{equation}\label{formal form}
f(z)=x+x^{m}+\sum_{j=m+1}^{\infty}a_{j}x^{j}.
\end{equation}
Then, there exists $\mu \in K$ and a formal power series $H(x)=x+\cdots $ so that $H\circ f\circ H^{-1}(x)=x+x^{m}+\mu x^{2m-1}$.
\end{prop}
\begin{proof}
We consider polynomials $h_{n}(x)=x+c_{n}x^{n}$, and define inductively $H_{2}(x)=h_{2}(x)$, and $H_{n}(x)=h_{n}\circ H_{n-1}(x)$ for $n >2$. Let us consider the effect of $H_{n}$ on $f$ for small $n$. Let $g(x)=x+x^{m}+b_{m+1}x^{m+1}+\cdots $. We shall show that the coefficients $c_{n}$ can be chosen so that the functions $H_{n}\circ f\circ H_{n}^{-1}$ and $g$ agree up to some specified order. To that end, let us consider $H_{2}\circ f\circ H_{2}^{-1}$. This agrees with $g$ up to order $m+1$, provided that
\begin{equation}
c_{2}=\frac{a_{m+1}-b_{m+1}}{m-2}.
\end{equation}
Of course, this change of variable will generally have an effect on each of the higher order terms, but that is of no concern to us formally. Thus, from now on, we consider $\alpha _{n+m-1}$ to be the coefficient of the $(n+m-1)$-degree term in $H_{n-1}\circ f\circ H_{n-1}^{-1}$. With this stipulation, we obtain the general formula for the coefficients $c_{n}$:
\begin{equation}\label{equation formula}
c_{n}=\frac{\alpha _{n+m-1}-b_{n+m-1}}{m-n}.
\end{equation}
Note, however, that this process breaks down when $n=m$; that is, the $(2m-1)$-degree term cannot be altered by these means. Thus, we see the invariant $\mu $ that appears; it is simply the coefficient $\alpha _{2m-1}$ in the expansion of the function $H_{m-1}\circ f\circ H_{m-1}^{-1}$. 

Finally, the formal map $H$ is defined to be $H=\lim_{n\rightarrow \infty}H_{n}$. Since the $n$th coefficient of $H_{n}$ is unchanged for all $H_{l}$ with $l>n$, we see that each coefficient in the formal series $H$ depends algebraically on a finite number of terms, and thus is well-defined. This completes the proof.
\end{proof}

\noindent {\it Remark}: Since the coefficient $c_m$ has no effect on the process outlined above, it can be considered a ``free term''. In what follows, we will take $c_m=0$. One consequence of the formal classification is that any mapping $f$ of the form (\ref{formal form}) can be taken to the form $\tilde f(x)=x+x^{m}+\mu x^{2m-1}+O(x^{2m})$ by a polynomial change of variable of degree $m-1$, and moreover, the proof shows that this change of variable is unique, if chosen so that it is tangent to the identity. Therefore, in much of what follows, we will assume that $f(x)=x+x^{m}+\mu x^{2m-1}+\cdots $, and therefore that $H(x)=x+A_{m+1}x^{m+1}+\cdots $. Finally, in the case that $K={\mathbf C}$, there is a formula relating the coefficient $\mu $ to a certain integral, but it is not of interest to us. Nonetheless, it is worth mentioning that in the case $m=2$, we have that $\mu =\frac{ a_{3}}{a_{2}^{2}}$. 

\subsection{Miscellaneous Notation}

We will often need to study the process of raising a power series to a given exponent (multiplicatively).
Consider, for example the problem of raising the power series
\[ f(x)=\sum_{i=0}^\infty a_{i}x^{i}\]
to the power $\ell$.

Then $f(x)^\ell$ will be a sum of terms of the form
\[   a_{i_1}a_{i_2} \cdots a_{i_\ell} x^{i_1+i_2+\cdots+i_\ell}, \]
where $i_1, i_2, \ldots, i_\ell$ is a (finite) sequence of positive integers, not necessarily distinct.

\begin{define} \label{length} Given a finite sequence $\ul{i}=(i_1,\ldots,i_\ell)$, write $|\ul{i}|=i_1+\cdots+i_\ell$. 
Also write $\ell(\ul{i})=\ell$, the ``length'' of $\ul{i}$.
\end{define}

We also adopt the following notation.
\begin{define} \label{bracket} Given a power series $f$, we write $[f]_{n}$ to be the coefficient of $x^{n}$. \end{define}

Thus, any power series $f$ may be written as $f(x)=\sum_n [f]_nx^n$.

The following lemma will be useful later.  Its proof is immediate.

\begin{lemma} \label{bracket lemma} Let $\eta(x)=\alpha_1x+ \alpha_d x^d + \alpha_{d+1}x^{d+1}+ \cdots \in K[[x]]$, and $j,T \in \N$.  Then if $j \neq T$ and $T<j+d-1$, then $[\eta(x)^j]_T=0$. \end{lemma} 

\section{The Integral Case}\setcounter{equation}{0}\label{section:s3}
We consider the case of mappings tangent to the identity. Before proving Theorem \ref{thm1.1}, we consider the case where the coefficients are in $\Delta$. Interesting in its own right, the proof of this case (and specifically, the proof of Proposition \ref{special case of theorem}) provides techniques which are crucial in the handling of the general case, proven in Section \ref{section:s4} (and in particular, in the proof of Proposition \ref{new c_{n} estimate}).

Let us first review the complex setting. As mentioned in Section \ref{section:s1}, there is a full holomorphic classification of mappings tangent to the identity in ${\mathbf C}$. However, the invariants are impractical for precise examples of conjugacy; in particular, there is no efficient way to determine if two formally equivalent mappings are tangent to the identity. However, some partial results do exist. We consider one such result. Let $f_{1,m}(z)=z+z^{2}$. While the full equivalence class of $f_{1,m}$ remains unknown, Ahern and Rosay \cite{A-R} have shown that the only germs of {\it entire} mappings which are holomorphically conjugate to $f_{0,m}$ are the mappings $f_{a}(x)=x+ax^{m}$, where $a\in \mathbf{C}$. Note that the conjugating map $h$ sending $f_{a,m}$ to $f_{1,m}$ takes the form $h(z)=az$. In particular, if we restrict the set of conjugating maps $h$ satisfying $h\circ g\circ h^{-1}=f_{1,m}$ to those which are tangent to the identity, then there are no entire mappings which are holomorphically equivalent to $f$. This is  in stark contrast to our own Theorem \ref{thm1.1}. Indeed, our theorem shows immediately that any polynomial of the form $p(x)=x+x^{m}+a_{2m}x^{2m}+\cdots +a_{k}x^{k}$ will be analytically equivalent to $f_{1,m}$.

For this section, we will consider series $f$ with coefficients in $\Delta $; write $f(x)=x+x^{m}+\cdots $. Such a mapping has as its normal form $f_{0,\mu }(x)=x+x^{m}+\mu x^{2m-1}$ (and since $\mu $ is determined by $f$, we will drop it for convenience and refer to $f_{0,\mu }$ as $f_{0}$). 

We need to study how the power series $H$ in the proof of Proposition \ref{formal equivalence} combines the coefficients $c_j$ of the polynomials $h_j$.  The following lemma, which is purely algebraic, determines which products may occur in a given degree. 

Let $R$ be a ring.  Let $\ul{c}=c_2,c_3, \ldots$ be a sequence of indeterminates.   
Write $\mathcal{A}=R[c_2, \ldots]$ for the polynomial ring in the variables $\{c_i \}$.
Suppose that $\ul{i}=(i_1,\ldots, i_\ell)$ is a finite sequence of natural numbers (not necessarily distinct).  
We define $\ell(\ul{i})=\ell$ and $|\ul{i}|=i_1+\ldots+ i_\ell$ as in Definition \ref{length}.
Write $c_{\ul{i}}$ for the monomial $c_{i_1} \cdots c_{i_\ell} \in \mathcal{A}$; its degree is $|\ul{i}|$.

Then a typical element of $\mathcal{A}$ may be written as $a(\ul{c})=\sum_{\ul{i}}\alpha_{\ul{i}} c_{\ul{i}}$, with 
$\alpha_{\ul{i}} \in R $.  

\begin{lemma}\label{lemma 3.1} 
For $j \geq 2$ let $h_j(x)=x+c_jx^j \in \mathcal{A}[x]$, and $H_j=h_j \circ h_{j-1} \circ \cdots \circ h_2 \in \mathcal{A}[x]$.  
Write $H_j(x)=x+\sum_n A_n^j(\ul{c})x^n$, with $A_n^j(\ul{c})$.  Suppose for a given $n \geq 2$, $A_n^j(\ul{c})=\sum_{\ul{i}}\alpha_{\ul{i}}^j c_{\ul{i}}$, with 
$\alpha_{\ul{i}}^j$ nonzero integers.  Then $n=|\ul{i}|-\ell(i)+1$. 
\end{lemma}

\begin{proof} 
We induct on $j$.  The statement is clear if $j=2$.
Given   $\ul{i}$, let $n(\ul{i})=|\ul{i}|-\ell(\ul{i})+1$.
Then $H_{j+1}(x)=H_j(x)+c_{j+1}(H_j(x))^{j+1}$.  The terms of the first part, $H_j(x)$, satisfy the proposition by induction.
By our inductive hypothesis, the second part is a sum of monomials of the form 
\[ c_{j+1} (\alpha_{\ul{i}_1}^j c_{\ul{i}_1}) \cdots (\alpha_{\ul{i}_L}^j c_{\ul{i}_L}) x^{\sum_{k=1}^L n(\ul{i}_k)} \cdot 
x^{(j+1)-L} .  \]

Therefore the exponent of $x$ in this monomial is given by
\[\left(\sum_{k=1}^L n(\ul{i}_k) \right) +(j+1)-L=\sum_{k=1}^L (|\ul{i}_k|-\ell(\ul{i}_k)+1)+(j+1)-L=  \sum_{k=1}^L (|\ul{i}_k|-\ell(\ul{i}_k))+(j+1) .            \]

On the other hand, write $\ul{i}'$ for the new sequence formed by concatenating $j+1, \ul{i}_1,\cdots$, and $\ul{i}_L$.
We have 
\[n(\ul{i}')=\left((j+1)+\sum_{k=1}^L |\ul{i}_k| \right)- \left(1+\sum_{k=1}^L \ell(\ul{i}_k)\right)+1,\]
which is equal to the previous expression.
\end{proof}

Now suppose $K$ is a non-archimedean field, with norm $|\cdot |$, and let $\Delta $ again denote those elements $x\in K$ with $|x|\leq 1$. Given a map $f\in K[|x|]$, we will write that $f\in \Delta [|x|]$ if all of the coefficients of $f$ lie in $\Delta $.

\begin{lemma}\label{H_{n}-coefficient estimate} 
Fix a positive integer $m$.  Let $c_j \in K$ for $j=m+1,\ldots ,n$, with $(j-m)!c_j \in \Delta $.  Let $h_j(x)=x+c_jx^j \in K[x]$, and $H_j=h_j \circ h_{j-1} \circ \cdots h_{m+1} \in K[x]$.  Write $H_j(x)=x+\sum_n A_nx^n$.  Then $(n-m)!A_n \in \Delta $. 
\end{lemma}

\begin{proof} 
By Lemma \ref{lemma 3.1}, we know that $H_j(x)=x+\sum_n A_n^j(\ul{c})x^n$, whose $n$th term is 
$A_n^j(\ul{c})=\sum_{\ul{i}}\alpha_{\ul{i}}^j c_{\ul{i}}$.  
The coefficients $\alpha_{\ul{i}}^j$ are integers, which will be nonzero only when $n=|\ul{i}|-\ell(\ul{i})+1$.  
So for the $n$th term we need only consider products of the form $c_{i_1} \cdots c_{i_\ell}$, with
\[ n=(i_1+\cdots + i_\ell)-\ell+1. \]
By hypothesis we have $(i_1-m)!\cdot (i_\ell-m)! c_{i_1} \cdots c_{i_\ell} \in \Delta $.  We know that the multinomial coefficient

\[\binom{i_1+i_2+ \cdots +i_\ell - \ell m}{(i_1-m),(i_2-m),\ldots, (i_\ell-m)}  =\frac{ ( |\ul{i}| - \ell m)  !}{  (i_1-m)!(i_2-m)!\cdots (i_\ell-m)!}      \]
is an integer.  It is therefore enough to prove that $n-m \geq |\ul{i}| - \ell m$.
By the equation for $n$ this reduces to showing that 
\[ |\ul{i}|-\ell+1 \geq |\ul{i}| - \ell m, \]
which is true since $\ell \geq 1$.

\end{proof}

We can even bound coefficients of powers of the $H_j(x)$ with similar methods.
\begin{lemma}\label{powers}
Let $N$ be a natural number, and write $H_j(x)^N=x^N+\sum_{n \geq m+1} b_nx^n$.  Then $(n-m)!b_n \in \Delta $.
\end{lemma}

\begin{proof} Similar to Lemma \ref{H_{n}-coefficient estimate}
\end{proof}

To recall, let $f\in \Delta [[x]]$ be a power series of the form (\ref{formal form}). By the remark following Proposition \ref{formal equivalence}, there is a polynomial change of variable which is tangent to the identity conjugating any such series to $f(x)=x+x^{m}+\mu x^{2m-1}+O(x^{2m}) .$ Thus, we can assume that $f$ takes the form 
\begin{equation}\label{reduced formal form}
f(x)=x+x^{m}+\mu x^{2m-1}+\sum_{n=2m}^{\infty }a_{n}x^{n} , 
\end{equation}
so that $f$ is formally equivalent to $f_{0}(x)=x+x^{m}+\mu x^{2m-1}$. Then, via Proposition \ref{formal equivalence} there is a formal series $H(x)=x+c_{m+1}x^{m+1}+\ldots $ conjugating $f$ with $f_{0,m}$, where $c_{k}$ is given by Equation (\ref{equation formula}) for all $k$ (and note that the series is unique, since we have chosen $c_{m}=0$). We show that this series converges in some neighborhood of $0\in K$. 

\begin{prop}\label{special case of theorem}
Let $f$ be an analytic mapping of the form (\ref{reduced formal form}), with coefficients $\mu \in \Delta $, $a_{n}\in \Delta $ for $n\geq 2m$. Let $h_{n}$, $H_{n}$, and $c_{n}$ be defined as in Proposition \ref{formal equivalence}. Then, $(n-m)!c_{n}\in \Delta $ for all $n\geq m+1$.
\end{prop}

\begin{proof}
We apply induction. For $n=m+1$, this follows immediately from (\ref{equation formula}) (and for $n=2,3,\ldots ,m$, we have chosen $c_{n}=0$). We now assume the statement is true for $n>m+1$, and prove it for $n+1$. Consider the equation $H_{n+1}\circ f=f_{0,m}\circ H_{n+1}$. By the definition of $H_{n+1}$, we can write
\begin{equation}\label{finite order equation}
\begin{split}
(H_{n}\circ f) +c_{n+1}(H_{n}\circ f)^{n+1}=\ &(H_{n}+c_{n+1}H_{n}^{n+1})+(H_{n}+c_{n+1}H_{n}^{n+1})^{m}\\
& +\mu (H_{n}+c_{n+1}H_{n}^{n+1})^{2m-1},\\
\end{split}
\end{equation}
(and note again that the powers are in fact multiplicative powers and not compositional powers). We wish to show that $c_{n+1}$ satisfies the inductive hypothesis. To this end, let us compute the coefficient of $x^{n+m}$ on both sides of (\ref{finite order equation}), noting that the two sides of (\ref{finite order equation}) agree up to order $O(x^{n+m+1})$, via Proposition \ref{formal equivalence}. We write Equation (\ref{finite order equation}) in the form
\begin{equation}\label{finite order equation 2}
\begin{split}
c_{n+1}(H_{n+1}\circ f)^{n+1}=\ &(H_{n}+c_{n+1}H_{n}^{n+1})+(H_{n}+c_{n+1}H_{n}^{n+1})^{m}\\
& +(H_{n}+c_{n+1}H_{n}^{n+1})^{2m-1}-(H_{n}\circ f),\\
\end{split}
\end{equation}
The $(n+m)$-degree term on the left-hand side of Equation (\ref{finite order equation 2}) is easy to compute - it is simply $(n+1)c_{n+1}x^{n+m}$. We compute the contribution from the right-hand side in parts.  Using Lemma \ref{bracket lemma} we see that the only terms which can contribute a coefficient of order $n+m$ are the maps $H_{n}$, $H_{n}^{m}$, $\mu H_{n}^{2m-1}$, $mc_{n+1}H_{n}^{n+m}$, and $-(H_{n}\circ f)$. Moreover, the contribution from $mc_{n+1}H_{n}^{n+m}$ is simply $mc_{n+1}x^{n+m}$. Thus, we can write
\begin{equation}\label{finite order equation 3}
(n-m+1)c_{n+1}=[H_{n}+H_{n}^{m}+H_{n}^{2m-1}-H_{n}\circ f]_{n+m}.
\end{equation}
 Thus, to complete the induction, we must show that the right-hand side of Equation \ref{finite order equation 3} belongs to $\frac{1}{(n-m)!}\Delta $.

 We write 
\begin{equation}
\begin{split}
& H_{n}(x)=x+\sum_{j=m+1}^{n+m}A_{j}x^{j} +O(x^{n+m+1}),\\ 
&H_{n}\circ f(x)=f(x)+\sum_{j=m+1}^{n+m}A_{j}(f(x))^{j}+O(x^{n+m+1}).\\
\end{split}
\end{equation}

By Lemma \ref{H_{n}-coefficient estimate}, $(j-m)!A_{j}\in \Delta $.

We first study the integrality $[H_n \circ f-H_n]_{n+m}$.  Since $x,f(x) \in \Delta [[x]]$, by the above formulas this reduces to the integrality of the $(n+m)$-coefficient of 
\begin{equation} \label{ag}
 \sum_{j=m+1}^{n+m} A_{j}(f(x)^j-x^j). 
\end{equation}
Let $g_j(x)=f(x)^j-x^j$; note that $g_j \in O(x^{m+j-1}) \cap \Delta [[x]]$.
We may discard most of the terms $A_j g_j$ in our computation: For $j \leq n$, we know that $(j-m)!A_j \in \Delta $ and therefore
$(n-m)!A_j \in \Delta $.  Since $g \in \Delta [[x]]$, we obtain $(n-m)! [A_jg_j]_{n+m} \in \Delta $.   For $j \geq n+2$, we know that $g_j \in O(x^{n+m+1})$, and therefore $[A_jg_j]_{n+m}=0$.

What happens for $j=n+1$?  We compute that 

\begin{equation} \label{miracle}
 [A_{n+1}(f(x)^{n+1}-x^{n+1})]_{n+m}=(n+1)A_{n+1}. 
\end{equation}

This seems unfortunate, but  in fact we will show that $(n-m)! [A_{n+1}g_{n+1} - H_n^m]_{n+m} \in \Delta $.


It is time to examine  $[H_n(x)^m]_{n+m}$.

Expanding the $m$th power of $H_n(x)$ and subtracting $x^m$ gives a sum of terms of the form
\[ x^{i_0}(A_{s_1}x^{s_1})^{i_1} \cdots (A_{s_T}x^{s_T})^{i_T}, \]
for integers $i_0,\ldots,i_T$ with $i_0 \neq m$ satisfying 
\begin{equation}\label{EQ1}
i_0+ \cdots + i_T=m.
\end{equation}

We are interested in those terms of degree $n+m$, which means that these integers also satisfy
\begin{equation}\label{EQ2}
i_0+i_1s_1+ \cdots + i_Ts_T=n+m
\end{equation}
As in the proof of Proposition \ref{powers}, if
\[ i_1(s_1-m)+\cdots+i_T(s_T-m) \leq n-m, \]
then the product of $(n-m)!$ with this coefficient is integral.
In view of (\ref{EQ2}), this will be true exactly when 
\begin{equation} \label{who cares}
 2m \leq i_0+(i_1+\cdots+i_T)m. 
\end{equation}

This inequality certainly holds if $i_1+ \cdots+i_T \geq 2$.  Suppose this is not the case.
Then $T=1$ and $i_1=1$.  Moreover (\ref{EQ1}) tells us that $i_0=m-1$, and then (\ref{EQ2}) tells us that $s_1=n+1$.
 
The conclusion of the above reasoning is that we have reduced to the sum of the terms $x^{m-1}A_{n+1}x^{n+1}=A_{n+1}x^{n+m}$ in the $m$th power expansion of $H_n(x)$.  In fact, there are $m$ such terms in the expansion, giving the term $mA_{n+1}x^{n+m}$.

Miraculously this fixes the ornery term from $[H_n \circ f-H_n]_{n+m}$.
Combining this with the leftover term (\ref{miracle}) gives the last piece $-(n-m+1)A_{n+1}$ of $[H_n+H_n^m-H_n \circ f]_{n+m}$, and this extra coefficient is exactly what we need since $(n-m)!(n-m+1)A_{n+1} \in \Delta $.

Finally, the analysis of $[\mu H_{n}^{2m-1}]_{n+m}$ is similar to, but easier than, the analysis of $[H_{n}^{m}]_{n+m}$. The only difference is that (\ref{EQ1}) now becomes
\begin{equation}\label{new EQ1}
i_{0}+\cdots +i_{T}=2m-1 .
\end{equation}
By (\ref{who cares}), if $i_{1}+\cdots +i_{T}=1$, then $i_{1}=1$, $i_{0}=2m-1$, and so $s_{1}=n-m+2$. Since $m\geq 2$, we see that $s_{1}\leq n$, and thus the estimate is satisfied.

\end{proof}

\begin{cor} The power series $H(x)=\lim_{n\rightarrow \infty}H_{n}(x) \in \OK$. \end{cor}
\begin{proof}
A final application of Lemma \ref{H_{n}-coefficient estimate} ensures that $(j-m)![H]_{j}\in \Delta $ for all $j\geq m+1$. This series will therefore converge in some disc centered at $0$, since $|n!|^{-1}\leq \alpha ^{-n}$. Moreover, since $H$ is tangent to the identity, this series is invertible near $0$. 
\end{proof}

\section{Vector Fields and Flows in $K$}\setcounter{equation}{0}\label{section:s2.5}

This section is devoted to the formal and analytic theory of flows and vector fields in $K$. The results here may be known, but we nonetheless give a complete proof of Theorem \ref{thm1.2}, and list some other facts about vector fields and flows which will be useful to us later. Much of the formal theory cited can be found in \cite{Arnold}.

We recall that a flow parametrized by an abelian group $S$ is a family of mappings $\{ \varphi ^{t}\}_{t\in S} $ satisfying the following conditions:
\begin{enumerate}
\item
$\varphi ^{0} ={\text Id}$,
\item
$\varphi ^{t_1}\circ \varphi ^{t_2}=\varphi ^{t_1+t_2}$.
\end{enumerate}
In particular, the flow is itself a commutative group under composition (for us $S$ will be the additive group of $K$). Given a formal vector field $V(x)=\sum _{n=2}^{\infty }v_{n}x^{n} \frac{\partial}{\partial x} $ with coefficients in $K$, it generates a formal flow $\{ T_{V}^{t}\}$ defined by $T^{0}_{V}(x)=x$, and $T_{V}^{t}(x)=\sum _{n=1}^{\infty }a_{n}(t)x^{n}$ defined by the set of differential equations
\begin{equation} \label{flow}
\frac{dT^{t}_{V}}{dt}(x)=V(T^{t}_{V}(x)).
\end{equation}
The differential equations are purely formal, and may be solved via (formal) anti-differentiation. However, we show  that in the non-archimedean setting, this formal operation gives rise to convergent power series in some small neighborhood of the origin. 

We first consider Theorem \ref{thm1.2} for a vector field with coefficients $a_{n}\in \Delta$. Later we deduce the general case from this.

Let 
\begin{equation}\label{vector field form}
V(x)=\sum_{n=2}^{\infty }v_{n}x^{n} \frac{\partial}{\partial x}
\end{equation}
be an analytic vector field whose coefficients satisfy $v_{n}\in \Delta $. We have the following proposition:

\begin{prop}\label{coefficients of flow}
Let $T^{t}_{V}(x)=T_{V}(t,x)$ be the formal flow of a vector field $V$ of the form (\ref{vector field form}) with coefficients $v_{n}\in \Delta $. Then, this flow is analytic in $x$ for all $t$.
\end{prop}

\begin{proof}
We write 
\begin{equation*}
T^{t}_{V}(x)=\sum _{n=1}^{\infty }a_{n}(t)x^{n}.
\end{equation*}
As before, write $a_{\ul{i}}=a_{i_{1}}\cdots a_{i_{\ell }}$ (note that this is a function in $t$), and write $|\ul{i}|=i_{1}+\cdots +i_\ell$. Since $v_{n}\in \Delta $ for all $n\geq 2$, Equation (\ref{flow}) implies that 
\begin{equation*}
\frac{da_{n}}{dt}=\dot{a}_{n}=\sum _{|\ul{i}|=n}\alpha _{n}a_{\ul{i}},
\end{equation*}
where $\alpha _{n}\in \Delta $ for all $n$.  Note that $a_{1}(t)\equiv 1$, by the definition of the flow.  An easy induction shows that deg$(a_{n})\leq n-1$ for all $n$. 
We now estimate the growth of the coefficients $a_{n}(t)$. Write $[\dot{a}_{n}(t)]_{m}$ for the coefficient of $t^{m}$ in the expansion of $\dot{a}_{n}$.  We will prove the following two-part hypothesis by induction: both that $m![\dot{a}_{n}]_{m}\in \Delta $ and $n!a_{n}(t)\in \Delta [t]$. Since $a_{1}(t)=1$ for all $t$, the base case is trivial. We first consider
\begin{equation*}
[\dot{a}_{n+1}]_{m}=\sum _{\ul{i}=n+1}\alpha _{n}[a_{\ul{i}}]_{m}.
\end{equation*}
where $\alpha_{\ul{i}}$ is simply an integer. Note that $[a_{\ul{i}}]_{m}$ is a product
\begin{equation*}
[a_{\ul{i}}]_{m}=[a_{i_{1}}\cdots a_{i_{l}}]_{m}=[a_{i_{1}}]_{m_{1}}\cdots [a_{i_{\ell }}]_{m_{l}},
\end{equation*}
where $m_{1}+\cdots +m_{\ell }=m$. Since deg$(a_{m_{j}})<m_{j}-1$ for all $j$, we may apply the induction hypothesis, and so  $m_{1}!\cdots m_{\ell }![a_{\ul{i}}]_{m}$ is integral. Therefore, $m![\dot{a}_{n+1}]_{m}$ is integral. Further, since 
\begin{equation*}
\dot{a}_{n+1}(t)=\sum_{m=1}^{n-1}[\dot{a}_{n+1}]_{m}t^{m},
\end{equation*}
we antidifferentiate to conclude that $(n+1)!a_{n+1}(t)\in \Delta [t]$.  This completes the two-part induction. 

We are now able to write 
\[
T^{t}(x)=x+\sum_{n=2}^{\infty }\frac{b_{n}(t)}{n!}x^{n},
\]
where $b_{n}(t)\in \Delta [t]$ has degree less than or equal to $n-1$. Note also that $|b_{n}(t)|\leq 1$ if $|t|\leq 1$, while $|b_{n}(t)|\leq |t|^{n-1}$ if $|t|>1$. Thus, for any $t\in K$ one can choose $r>0$ so that if $x\in B_{r}(0)$, then the series $T^{t}_{V}$ converges.
\end{proof}

\noindent
{\it Proof of Theorem \ref{thm1.2}}. We consider an arbitrary vector field $V$ of the form (\ref{vector field form}). Let $r\in K$ and consider the change of variables $L_{r}(x)=rx$. The resulting vector field $\widetilde V$ has the form
\[
\widetilde V(x)=\frac{1}{r}\sum_{n=2}^{\infty }v_{n}(rx)^{n}\frac{\partial }{\partial x}.
\]
Since $V$ is analytic, one can choose $r$ so that $|r^{n-1}v_{n}|\leq 1$ for all $n$. Thus, by Proposition \ref{coefficients of flow}, the formal flow $T_{\widetilde{V}}^{t}(x)$ associated with $\widetilde V$ is in fact analytic, and for all values of time $t$, we have $L_{r}\circ g_{\widetilde{V}}^{t}\circ L_{r}^{-1}=g_{V}^{t}$, where $g^{t}_{V}$ (resp. $g^{t}_{\widetilde {V}}$) is the time-$t$ map of the flow $T^{t}_{V}$ (resp. $T^{t}_{\widetilde{V}}$). Thus, the formal flow of $V$ is also analytic, completing the proof.

We shall make use of the following lemma in Section \ref{section:s4}. The result is purely formal, and thus the elementary proof in ${\mathbf C}$ can be modified for fields of characteristic $0$.

\begin{lemma}\label{flows and centralizers}
Let $V(x)$ be a vector field which generates a formal flow $T^{t}_{V}(x)$. Fix $t=t_{0}$, and write $f(x)=T^{t_{0}}_{V}(x)$. Then, if $g(x)=x+\cdots $ is any formal map satisfying $g\circ f=f\circ g$, then there is a $t_{1}$ so that $g(x)=T^{t_{1}}_{V}(x)$; i.e. $g$ is in the flow of $V(x)$.
\end{lemma}

Note that in our setting, this implies that any formal map $g$ centralizing the time-$t$ map of an analytic vector field $V$ must itself be analytic.

\section{A proof of Theorem \ref{thm1.1} and remarks on centralizers and root extraction}\setcounter{equation}{0}\label{section:s4}
We now prove Theorem \ref{thm1.1}.  
Let $f(x)=x+x^m+ \mu x^{2m-1}+O(x^{2m}\in \OK$ as in (\ref{reduced formal form}).  Since the radius of convergence $\rho$ is positive, the sequence $\{ \frac{1}{\sqrt[n]{ |a_n| }} \}$ is bounded below by some $\eps >0$.  Pick $q \in K$ with $|q| \leq \eps$.  Then 
$b_n=a_n q^n \in \Delta$ for all $n$.

Thus, we are reduced to the study of series $f$ of the form
\begin{equation}\label{new reduced form}
f(x)=x+x^{m}+\sum_{n=2m-1}^{\infty }\frac{b_{n}}{q^{n}}x^{n},
\end{equation}
where $b_{n}\in \Delta $. The idea here will be to estimate the decay of the denominators in the coefficients $c_{n}$ and $A_{n}$ in our formal conjugating maps. 

We introduce a ``jump'' function which governs the growth of the coefficients of $H_n$. Fix $m \geq 2 $. Let us define, for $n\in {\mathbf N}$ with $n\geq m+1$ 
\begin{equation}\label{sigma}
\sigma _{m}(n)=\left( \frac{2m-1}{m-1}\right) (n-(m+1))+\varepsilon _{m}(n),
\end{equation}
where
\begin{equation}\label{epsilon}
\varepsilon _{m}(n)=2m+\left( 1-\frac{2m-1}{m-1}\right) r, \text{ for } (n-2)\equiv r{\text{mod}}(m-1), 0\leq r\leq m-2.
\end{equation}

This function measures the growth of the power of $q$ in the denominator of the coefficients of the formal conjugating maps. In particular, if $m=2$, the function $\sigma _{2}$ is simply linear, with constant slope $3$. If $m>2$, the function is more complicated: it has a constant slope of $1$ on the interval $[m+1,2m-1]$, but as $n$ moves from $2m-1$ to $2m$, it  ``jumps'' by $m+1$.
In fact, the function $\sigma_m$ continues this behavior over each subsequent interval of length $m-1$.
If $m>2$, this mimics closely the growth of the denominators in our conjugating map $H$.  When $m=2$, the denominator growth is somewhat smaller; nonetheless our estimates handle this case. The following lemma describes the behavior of the functions $\sigma_{m}$ and $\varepsilon_{m}$.

\begin{lemma}\label{properties of sigma}
Let the functions $\sigma _{m}$ and $\varepsilon _{m}$ be given by Equations (\ref{sigma}) and (\ref{epsilon}), resp., and let $n\geq m+1$. Then the following hold.
\begin{enumerate}
\item
$\varepsilon _{m}$ is a periodic function of period $m-1$ which is decreasing on $m+1\leq n\leq 2m-1$. Moreover, 
\begin{equation}
\frac{m^2}{m+1} \leq \varepsilon _{m}(n)\leq 2m.
\end{equation}
\item
$\sigma _{m}(n)$ is a strictly increasing, integer-valued function of $n$,  $\sigma _{m}(n+(m-1))=\sigma _{m}(n)+(2m-1)$, and 
\begin{equation}
\left(\frac{2m-1}{m-1} \right) n-\left(m+2+\frac{1}{m-1} \right) \leq \sigma _{m}(n)\leq \left(\frac{2m-1}{m-1} \right) n - \left(\frac{3m-1}{m-1} \right).
\end{equation}

\item
If $a,b\in {\mathbf N}$, and $b-a\geq m-1$, then $\sigma_{m}(b)-\sigma_{m}(a)\geq (b-a)+m$.
\item
Let $\underline{i}=(i_{1},\cdots ,i_{l})$ be an $\ell$-tuple of positive integers and let $n=|\underline{i}|-\ell+1$. Then,
\begin{equation}
\sum_{j=1}^{\ell}\sigma _{m}(i_{j})\leq \sigma _{m}(n).
\end{equation}
\end{enumerate}
\end{lemma}
\begin{proof}
The first two statements are elementary. If $b=a+(m-1)$, then part (ii) of the lemma gives that
\begin{equation*}
\sigma_m(b)=\sigma_{m}(a)+(2m-1)=\sigma_{m}(a)+(b-a)+m.
\end{equation*}
Since  for any $n\in {\mathbf N}$, we have $\sigma_{m}(n+1)-\sigma_{m}(n)\geq 1$, statement (iii) follows. For the last statement, we note that, using the definition of $\sigma_m$ and the formula for $n$, this reduces to proving that 
\[ \sum_j \eps_m(i_j) \leq m\left( \frac{2m-1}{m-1} \right) (\ell-1)+ \eps_m(n). \]
Suppose that $n-2$ has remainder $r$, and $i_j-2$ has remainder $r_j \mod (m-1)$.
Then using the definition of $\eps_m$ this reduces to proving that
\[ r \leq \sum_j r_j +(\ell-1). \]
But we know $n-2$ is congruent mod $(m-1)$ to $ \left(\sum_j r_j+2\ell \right) -\ell -1$, whose remainder mod $(m-1)$ is no greater than than itself.
\end{proof}

Now, let $f$ be of the form (\ref{new reduced form}) with formal invariants $m$ and $\mu =b_{2m-1}/q^{(2m-1)}$. Associated to $m$, we have the function $\sigma_{m}$; we drop the $m$ for convenience. We now prove propositions similar to those of Section \ref{section:s3}.

\begin{prop}\label{new H_{n} estimate}
Fix a natural number $m\geq 2$, and let $c_{j}\in K$ for $j=m+1,\cdots ,n$, such that $(j-m)!q^{\sigma (j)}c_{j}\in \Delta$. Define $h_{j}=x+c_{j}x^{j}\in K[x]$, and write $H_{j}(x)=h_{j}\circ\cdots \circ h_{m+1}(x)=x+\sum_{n}A_{n}x^{n}$. Then for all $n$, $(n-m)!q^{\sigma (n)}A_{n}\in \Delta$.
\end{prop}
\begin{proof}
Recalling notation from Proposition  \ref{H_{n}-coefficient estimate} , we know that $H_j(x)=x+\sum_n A_n^j(\ul{c})x^n$, whose $n$th term is 
$A_n^j(\ul{c})=\sum_{\ul{i}}\alpha_{\ul{i}}^j c_{\ul{i}}$.  
The coefficients $\alpha_{\ul{i}}^j \in \Z$ will be nonzero only when $n=|\ul{i}|-\ell(\ul{i})+1$.  

By hypothesis we have $(i_1-m)! \cdots (i_\ell-m)! q^{\sigma (i_1)+\cdots \sigma (i_\ell)} c_{i_1} \cdots c_{i_\ell} \in \Delta$. 
As in the proof of Proposition \ref{H_{n}-coefficient estimate} we may replace the product of factorials with $(n-m)!$. The power 
$q^{\sigma(n)}$ is handled by Part (iv) of Lemma \ref{properties of sigma}.
\end{proof}

Similarly, for the coefficients $c_{n}$, we have the following:
\begin{prop}\label {new c_{n} estimate}
Let $f$ be an analytic mapping of the form (\ref{new reduced form}), where $b_{n}\in K$. Let $h_{n}$, $H_{n}$, and $c_{n}$ be defined as in Proposition \ref{formal equivalence}. Then, $(n-m)!q^{\sigma (n)}c_{n}\in \Delta$ for all $n\geq m+1$.
\end{prop}

The idea of the proof is similar to that of Proposition \ref{special case of theorem}, but more involved, since we now have to measure the growth of the denominators $q^k$. We begin by proving a series of lemmas. The first of these lemmas consider the growth of denominators in the coefficients of a power series $H_{n}^{m}$ as measured by the function $\sigma _{m}$.

\begin{lemma}\label{H_{n}^{m}}
Let 
\begin{equation}\label{H_{n}}
P_{n}(x)=x+\sum_{k=m+1}^{\infty}A_{k}x^{k}
\end{equation}
be a power series, and let $q\in K$. Suppose that $q^{\sigma (k)}A_{k}\in \Delta  $ for all $n\geq m+1$. Then $q^{\sigma (n+1)}[P_{n}^{m}]_{n+m}\in \Delta$. 
\end{lemma}
\begin{proof}
We write $[P_{n}^{m}]_{n+m}x^{n+m}$ as a sum of terms of the form
\begin{equation*}
x^{k_{0}}(A_{i_{1}}x^{i_{1}})\cdots (A_{i_{l}}x^{i_{\ell}}),
\end{equation*}
where 
\begin{equation*}
m=k_{0}+\ell
\end{equation*}
and
\begin{equation*}
n+m=k_{0}+i_{1}+\cdots +i_{\ell}.
\end{equation*}
The estimate in our hypothesis yields
\begin{equation*}
q^{\sigma (i_{1})+\cdots +\sigma (i_{\ell})}A_{i_{1}}\cdots A_{i_{\ell}}\in \Delta,
\end{equation*}
and so we must prove 
\begin{equation*}
\sum_{s=1}^{\ell}\sigma (i_{s})\leq \sigma (n+1).
\end{equation*}

However, if we write $|\ul{i}|=i_{1}+\cdots +i_{\ell}$, we see that $n+1=|\ul{i}|-\ell+1$.

 Thus, we may apply part (iv) of Lemma \ref{properties of sigma}.
\end{proof}

The next lemma is similar to Lemma \ref{H_{n}^{m}}, but we now consider coefficients of $H_{n}^{2m-1}$. However, the growth of the denominators here is actually a bit less than that in Lemma \ref{H_{n}^{m}}; as we shall see, this small improvement is crucial to the proof of Proposition \ref{new c_{n} estimate}.

\begin{lemma}\label{H_{n}^{2m-1}}
Let $P_{n}$ be as in (\ref{H_{n}}), and suppose that $A_{n}$ satisfies the estimates of Lemma \ref{H_{n}^{m}}. Then, we have $q^{\sigma (n+1)-(2m-1)}[P^{2m-1}_{n}]_{n+m}\in \Delta $.
\end{lemma}
\begin{proof}
Similar to the proof of the previous lemma, we can write $[P_{n}^{2m-1}]_{n+m}x^{n+m}$ as a sum of terms of the form
\begin{equation*}
x^{k_{0}}(A_{i_{1}}x^{i_{1}})\cdots (A_{i_{\ell}}x^{i_{\ell}}),
\end{equation*}
where 
\begin{equation*}
2m-1=k_{0}+\ell
\end{equation*}
and
\begin{equation*}
n+m=k_{0}+i_{1}+\cdots +i_{\ell}.
\end{equation*}
From the estimate in the hypothesis of Lemma \ref{H_{n}^{m}}, we have that 
\begin{equation*}
q^{\sigma (i_{1})+\cdots +\sigma (i_{\ell})}A_{i_{1}}\cdots A_{i_{\ell}}\in \Delta.
\end{equation*}
Thus, to prove the lemma, we must show that 
\begin{equation*}
\sum_{s=1}^{\ell}\sigma (i_{s})\leq \sigma (n+1)-(2m-1).
\end{equation*}
Again writing $|\ul{i}|=i_{1}+\cdots +i_{\ell}$, we have 
\begin{equation*}
|\ul{i}|=n+m-k_{0}=n+m-[(2m-1)-\ell].
\end{equation*}
Hence, we have 
\begin{equation*}
n+1=(|\ul{i}|-\ell+1)+(m-1),
\end{equation*}
and so, by part (ii) of Lemma \ref{properties of sigma}
\begin{equation*}
\sigma (n+1)=\sigma (|\ul{i}|-\ell+1)+(2m-1).
\end{equation*}
Hence, by part (iv) of Lemma \ref{properties of sigma}, the proof is complete.
\end{proof}

Our final lemma again estimates the denominators of coefficients - this time, we consider maps of the form $P_{n}-P_{n}\circ f$. In the proofs of the previous two lemmas, we really only needed to relate different values of $\sigma $. However, in the proof below, we must relate the growth of denominators which depend on $\sigma $ to those which do not depend on it - namely, the growth of the denominators of powers of $f$.

\begin{lemma}\label{H_{n}circ f}
Let $P_{n}$ be as in (\ref{H_{n}}), and suppose that $A_{n}$ satisfies the estimates of Lemma \ref{H_{n}^{m}}. Suppose that $f$ is of the form (\ref{new reduced form}). Then, we have $q^{\sigma (n+1)}[P_{n}-P_{n}\circ f]_{n+m}\in \Delta $.
\end{lemma}
\begin{proof}
It is enough to prove that for all $s$, 
\[ q^{\sigma(n+1)}[A_s(f(x))^s]_{n+m} \in \Delta. \]
By our estimate for $A_s$ we need only prove that
\[ q^{\sigma(n+1)-\sigma(s)}[(f(x))^s]_{n+m} \in \Delta. \]
Expanding the $s$ power of $f$ gives a sum of terms of the form
\begin{equation}
 x^{e_1}(x^m)^{e_m}\left( \frac{b_{2m-1}x^{2m-1}}{q}\right) ^{e_{2m-1}} \ldots \left( \frac{b_\ell x^\ell}{q}\right) ^{e_\ell}, 
\end{equation}
with 
 
\beq \label{veggirl}
 e_1 + e_m+ e_{2m-1}+\ldots + e_\ell=s. 
\end{equation} 

This term will have degree $n+m$ when
\beq \label{vg2}
 e_1 + me_m+(2m-1)e_{2m-1}+ \ldots+ \ell e_\ell=n+m.  
\end{equation}
We know that 
\[ q^{(2m-1)e_{2m-1}+ \ldots+ \ell e_\ell} \left[ x^{e_1}(x^m)^{e_m}\left( \frac{b_{2m-1}x^{2m-1}}{q}\right) ^{e_{2m-1}} \ldots \left( \frac{b_\ell x^\ell}{q}\right) ^{e_\ell} \right] _{n+m} \in \Delta, \] 
and therefore it is enough to prove that 

\beq \label{vg3}
 \sigma(n+1)-\sigma(s) \geq (2m-1)e_{2m-1}+ \ldots + \ell e_\ell.  
\end{equation}

By Lemma \ref{properties of sigma}, we have
\beq
\sigma(n+1)-\sigma(s) \geq \frac{2m-1}{m-1}(n+1-s)-\left( (m-1)+\frac{1}{m-1}\right) .
\end{equation}
Subtracting Equation \ref{veggirl} from Equation \ref{vg2} gives
\beq
n-s=(m-1)e_m+(2m-2)e_{2m-1}+\cdots+ (\ell-1)e_\ell-m. 
\end{equation}
Combining this with the above estimate, we see that the inequality (\ref{vg3})
will be true when
\beq 
(2m-1)e_m+(2m-1)e_{2m-1}+ \cdots + \left( (\ell-1)\left( \frac{2m-1}{m-1}\right) -\ell\right) e_\ell \geq (3m-2)-\frac{1}{m-1}. 
\end{equation}
In fact, all the  coefficients on the left hand side are greater or equal to $2m-1$.
We may therefore conclude the following:  If the inequality (\ref{vg3}) does {\it not} hold, then only one of the exponents 
$e_m, \ldots, e_\ell$ may be nonzero, and in fact must be equal to $1$.

Let us focus now on these simple terms; they must be of the form
\[ x^{s-1}b_{n+m-s+1} \left(\frac{x}{q} \right)^{n+m-s+1}, \]
with $n+m-s+1 \geq 2m-1$ (recall there is no denominator for the coefficient of $x^m$ in $f$).
Thus we only need to check that if $m-1 \leq n+1-s,$ then $\sigma(n+1)-\sigma(s) \geq (n+1-s)+m$. This follows immediately from statement (iii) of Lemma \ref{properties of sigma}.
\end{proof}

\noindent
{\it Proof of Proposition \ref{new c_{n} estimate}}.
We induct on $n$. For $n=m+1$, it is clear from the proof of Proposition \ref{formal equivalence} that $c_{m+1}$ satisfies the estimate, and we take $c_{2}=c_{3}=\cdots =c_{m}=0$. Thus, we assume that $c_{n}$ satisfies the estimate hypothesis, and we show that $c_{n+1}$ also satisfies it.

From Proposition \ref{new H_{n} estimate}, we can write $H_{n}$ in the form of (\ref{H_{n}}), where $q^{\sigma (n)}(n-m)!A_{n}\in \Delta _{p}$ for all $n\geq m+1$. In order to complete the proof of the proposition, we must show that $q^{\sigma (n+1)}(n-m+1)!c_{n+1}\in \Delta $.

Writing $H_{n+1}=H_{n}+c_{n+1}H_{n}^{n+1}$, the formal classification theorem shows that, up to order $O(x^{n+m+1})$, we must have
\begin{equation*}
\begin{split}
H_{n}\circ f+c_{n+1}(H_{n}\circ f)^{n+1}=
& H_{n}+c_{n+1}H^{n+1}+(H_{n}+c_{n+1}H_{n}^{n+1})^{m}\\
& +\frac{b_{2m-1}}{q^{2m-1}}(H_{n}+c_{n+1}H_{n}^{n+1})^{2m-1}.\\
\end{split}
\end{equation*}
We consider the $(n+m)$-degree coefficient of each side. Again, this reduces to the equation
\begin{equation}\label{new coefficient equation}
(n-m+1)c_{n+1}x^{n+m}=[H_{n}-H_{n}\circ f]_{n+m}+[H_{n}^{m}]_{n+m}+[H_{n}^{2m-1}]_{n+m}.
\end{equation}
We show that the terms on the right-hand side of (\ref{new coefficient equation}) each lie in $\frac{q^{-\sigma(n+1)}}{(n-m)!} \Delta $. The idea is to synthesize the lemmas proved above with the propositions proved in Section \ref{section:s3}. Let us first consider $[H_{n}^{2m-1}]_{n+m}$. We can write the $(n+m)$-degree term as a sum of terms of the form $x^{k_{0}}\prod_{s}A_{i_{s}}x^{i_{s}}$, and we know that $\prod_{s}(i_{s}-m)!q^{\sigma (i_{s})}A_{i_{s}}\in \Delta $. Since $i_{s}\leq n$ for all $s$, we may combine the techniques of Proposition \ref{special case of theorem} with those of Lemma \ref{H_{n}^{2m-1}} to conclude that $(n-m)!q^{\sigma (n+1)}[H_{n}^{2m-1}]_{n+m}\in \Delta $. Similarly, we consider the coefficient $[H_{n}^{m}]_{n+m}$. Again, write the $(n+m)$-degree term as a sum of terms of the form $x^{k_{0}}\prod_{s}A_{i_{s}}x^{i_{s}}$. As long as $i_{s}\leq n$, we can again apply the techniques of Proposition \ref{special case of theorem} and Lemma \ref{H_{n}^{m}} to show that any such term satisfies the necessary condition. However, while it is clear that $i_{s}\leq n+1$ for all $s$, it can occur that $k_{0}=m-1$ and $i_{1}=n+1$. The term associated with this takes the form $mA_{n+1}x^{n+m}$; note that $q^{\sigma (n+1)}(n-m+1)!A_{n+1}\in \Delta $.

Finally, we consider $[H_{n}-H_{n}\circ f]_{n+m}$.  It is easy to see that the $(n+m)$-degree term can be written as a sum of terms of the form $x^{k_{0}}\prod_{s}A_{i_{s}}x^{i_{s}}$, where $i_{s}\leq n+1$. Moreover, if $i_{s}=n+1$ for any $s$, then $s=1$, and the term takes the form $-(n+1)A_{n+1}x^{n+m}$. Thus, combining this with the leftover term from $[H_{n}^{m}]_{n+m}$, we obtain our ``miracle'' $-(n-m+1)A_{n+1}x^{n+m}$, which satisfies the necessary estimate. Finally, if $i_{s}\leq n$ for all $s$, then we can handle the estimate in the same manner as above, this time combining the techniques of Proposition \ref{special case of theorem} with Lemma \ref{H_{n}circ f}.

\qed
\\
\noindent
{\it Proof of Theorem \ref{thm1.1}.} We write $H_{n}(x)=(h_{n}\circ h_{n-1}\circ \cdots \circ h_{m+1})(x)$, where $h_{k}(x)=x+c_{k}x^{k}$. From Propositions \ref{new H_{n} estimate} and \ref{new c_{n} estimate}, we note that the conjugating map $H=\lim _{n\rightarrow \infty}H_{n}$ will have coefficients $A_{n}$ satisfying 
\begin{equation*}
(n-m)!q^{\sigma _{m}(n)}A_{n}\in \Delta,
\end{equation*}
where $q$ is chosen as in the beginning of Section \ref{section:s4}. By Proposition \ref{factorial estimate}, we have $|(n-m)!|^{-1}<\alpha ^{-n}$ for some real number $0<\alpha \leq 1$, and thus by choosing $q$ to satisfy also $0<|q|<\alpha $, we obtain that $|(n-m)!|^{-1}<|q|^{-n}$. From Lemma \ref{properties of sigma}, we have that $\sigma_{m}(n)\leq 3n$ for all $m\geq 2$, $n\geq m+1$. Thus, $A_{n}x^{n}$ will tend to $0$ if $|x|$ is sufficiently small, and hence our series converges.
\qed
\\

With a full analytic classification in place, we now settle the questions of centralizers and root extraction for a typical analytic map $f$ of the form (\ref{reduced formal form}).

Let us begin with centralizers. It is obvious that if $G$ is a group and $f$ and $g$ are conjugated by an element $h$, i.e. $h\circ f=g\circ h$, then if we write $Z(f)$ and $Z(g)$ for the centralizers of $f$ and $g$, respectively, then $h[Z(f)]h^{-1}=Z(g)$. Thus, given an analytic mapping $f$, we must find one map in the analytic equivalence class of $f$ for which the analytic centralizer is understood. A starting point would be a mapping for which the formal and analytic centralizer agree. 

Let $f$ be of the form (\ref{reduced formal form}) with formal invariants $m$ and $\mu $. We write 
\begin{equation*}
Z_{F}(f)=\{ g\in K[[x]]: g'(0)=1, g\circ f=f\circ g\},
\end{equation*}
\begin{equation*}
Z_{A}(f)=\{ g\in \OK :g'(0)=1, g\circ f=f\circ g\}.
\end{equation*}
Consider the vector field $V(z)=\frac{x^{m}}{1-(\mu -1) x^{m-1}}\frac{\partial }{\partial x}$. The time-one map $T^{1}_{V}$ of this vector field takes the form
\[ 
T^{1}_{V}(x)=x+x^{m}+\mu x^{2m-1}+\cdots. 
\]
Thus, $f$ and $T^{1}_{V}$ are analytically conjugate. Note also that $Z_{F}(T^{1}_{V})=Z_{A}(T^{1}_{V})=\{ T^{t}_{V}:t\in K\} $, by Lemma \ref{flows and centralizers} and the remarks following. Let $h$ be an analytic map tangent to the identity satisfying $h\circ f\circ h^{-1}=T^{1}_{V}$. We have the following:
\begin{cor}\label{centralizers}
$Z_{A}(f)=\{ h^{-1}\circ T^{t}_{V}\circ h \}_{t \in K}.$
\end{cor}
\begin{proof}
This follows immediately from Lemma \ref{flows and centralizers}.
\end{proof} 
Note also that $Z_{A}(f)=Z_{F}(f)$.

Root extraction is now a simple consequence. 
\begin{cor}\label{roots}
Let $f$ be of the form (\ref{reduced formal form}), and let $n\geq 1$ be a natural number. Then, there is a unique $g\in \OK$ tangent to the identity satisfying $g^{\circ n}=f$.
\end{cor}
\begin{proof}
Basic computation yields that there is a unique {\it formal} map $g$ tangent to the identity satisfying $g^{\circ n}=f$. Any such root necessarily belongs to the centralizer of $f$, and since the formal centralizer agrees with the analytic one, we conclude that $g$ is in fact analytic. Thus, $f$ admits analytic $n$th-root extraction of all orders. 
\end{proof}

\bibliographystyle{plain}

\begin{thebibliography}{99}
\bibitem{A-R} P. Ahern and J.-P. Rosay,
{\em Entire functions, in the classification of differentiable germs tangent to the identity, in one or two variables},
Trans. of the Amer. Math. Soc. {\bf 347} (1995), no. 2, 543-572.
\bibitem{Arnold} V. I. Arnold, 
{\bf Geometrical Methods in the Theory of Ordinary Differential Equations}, Springer-Verlag, New York, 1983.
\bibitem{Benedetto} R. Benedetto,
{\em Hyperbolic maps in $p$-adic dynamics},
Ergodic Th.  Dyn. Sys., {\bf 21} (2001), 1-11.
\bibitem{Camacho} C. Camacho,
{\em On the local structure of conformal mappings and holomorphic vector fields in ${\mathbf{C}}^{2}$}, 
Bull. Soc. Math. de France Ast\'{e}risque, (1978).
\bibitem{Ecalle} J. Ecalle
{\em Sur les functions r\'{e}surgentes}, 
I, II, Publ. Math. d'Orsay, Universit\'{e} de Paris-Sud, Orsay, 1981.
\bibitem{Fatou} P. Fatou,
{\em Sur les equation fonctionnelles},
Bull. Soc. Math. France {\bf 47} (1919), 161-271.
\bibitem{HC} Harish-Chandra,
{\em Harmonic Analysis on Reductive $p$-adic Groups},
Proc. of Symp. in Pure Math., {\bf XXVI} (1973), Amer. Math. Soc., Providence, R.I., 167-192.
\bibitem{Ilya} Y. S. Il'yashenko,
{\em Nonlinear Stokes Phenomena}.  
Adv. in Soviet Math., vol. 14, Amer. Math. Soc., Providence, RI, 1993.
\bibitem{Jenkins1} A. Jenkins,
{\em Further Reductions of Poincar\'{e}-Dulac Normal Forms},
to appear, Proc. of the Amer. Math. Soc.
\bibitem{Jenkins2} A. Jenkins,
{\em Holomorphic germs and smooth conjugacy in a punctured neighborhood of the origin},
Trans. of the Amer. Math. Soc. {\bf 360} (2008), no. 1, 331-346.
\bibitem{koenigs}G. Koenigs,
{\em Recherches sur les Integral de Certain Equation Fonctionelles}, 
Ann. Scient. Ec. Norm. Sup. {\bf 1} (1884), 1-41.
\bibitem{Malgrange} B. Malgrange
{\em Travaux d'Ecalle et de Martinet-Ramis sur les syst\` {e}mes dynamiques}, 
S\'{e}minar Bourbaki, vol. 1981/1982, Ast\'{e}risque, vol. 92-93, Soc. Math. France, Paris, 1982, pp. 59-73
\bibitem{Serre} J. P. Serre
{\bf Lie Algebras and Lie Groups: 1964 Lectures Given at Harvard University}, W.A.Benjamin Inc., New York, 1965.
\bibitem{Shch} A. A. Shcherbakov,
{\em Topological classification of germs of conformal mappings with identical linear part}, 
Vestnik Moskov. Univ. Ser. I Mat. Mekh., {\bf 1982}, no. 3, 52-57; English transl. in Moscow Univ. Math. Bull. {\bf 37} (1982).
\bibitem{Schikhof} W. H. Schikhof,
{\bf Ultrametric Calculus: An Introduction to $p$-adic Analysis},
Cambridge Studies in Advanced Mathematics, {\bf 4}, Cambridge University Press, Cambridge, 1984.
\bibitem{Voronin} S. M. Voronin,
{\em Analytic classification of germs of conformal maps $({\mathbb{C}},0)\rightarrow ({\mathbb{C}},0)$ with identical linear part}, Func. Anal. Appl. {\bf 15} (1981), 1-17. 
\end{thebibliography}

\end{document}